\documentclass[12pt,Йeqno]{article}
\usepackage{cmap}
\usepackage[koi8-r]{inputenc}
\usepackage{graphicx}
\usepackage{natbib}
\usepackage{amsfonts}
\usepackage[tbtags]{amsmath}
\topskip=\baselineskip%
\newcommand{\lr}{\left(}

\newcommand{\rr}{\right)}
\newcommand{\lf}{\left\{}
\newcommand{\rf}{\right\}}
\newcommand{\lv}{\left|}
\newcommand{\rv}{\right|}
\newcommand{\lp}{\left.}
\newcommand{\rp}{\right.}

\newcommand{\eel}{\end{lemma}}

\newcommand{\xin}{X_{i:n}}

\newcommand{\al}{\alpha}
\newcommand{\be}{\beta}

\title{{\Huge \sf Bounds for absolute moments of order statistics}
\author{\large Nadezhda V. Gribkova}
\date{
{\small \it
\centerline{St.Petersburg~State~University,~Mathematics~and~Mechanics~Faculty,}
\centerline{199034,~ Universitetskaya nab.~7/9, St.~Petersburg, Russia, n.gribkova@spbu.ru}
}}}
\begin{document}
\noindent{\emph{Exploring Stochastic Laws.~pp.~129--134.}}\\
\noindent{\emph{Eds.~A.V.\,Skorokhod and Yu.V.\,Borovskikh}}\\
\noindent{\copyright \emph{ VSP \,Utrecht, 1995}}\\

\bigskip

\bigskip
\noindent{\sf \LARGE Bounds for absolute moments of order statistics}

\bigskip

\bigskip
\noindent{NADEZHDA V. GRIBKOVA}

\bigskip
\noindent{\emph{Department of Applied Mathematics, Transport University,}}\\
\noindent{\emph{199031, St.\,Petersburg, Moskovsky Prospect 9, Russia}}

\bigskip

\bigskip
\noindent{Received February 23, 1995}

\medskip

\bigskip

\bigskip
\noindent{\bf Abstract.} {\small The bounds for absolute moments of order statistics are established.
Let $X_1,\dots ,X_n$ be independent identically distributed real-valued random variables and let $X_{1:n}\le \dots \le X_{n:n}$ be the corresponding order statistics. The absolute moments $\textbf{E}|\xin|^k$, $k>0$, are estimated via the absolute moment $\textbf{E}|X_1|^{\delta}$, $\delta>0$, for all $i$ such that $k\delta^{-1}\leq i \leq n-k\delta^{-1}+1$ with order $(n^2i^{-1}(n-i)^{-1})^{k\delta^{-1}}$ in $i$ and $n$. These estimates are able to be of some use as a~tool to argue in different probability limit theorems.}

\bigskip

\\[3mm] \noindent{\bf MSC:} \  62G30, 60E15.\\[3mm]

\noindent{Let $X_1,\dots ,X_n$ be independent identically distributed real-valued random variables with common distribution function $F$, and let $X_{1:n}\le \dots \le X_{n:n}$ denote  the corresponding order statistics.
Let $k$ and $\delta$ be arbitrary positive numbers. Put $\rho=k\delta^{-1}$ and set $g(u)=u(1-u)$.

\bigskip

\noindent{\bf Theorem 1.} {\it  For all $n\geq 2\rho+1$ and for all $i$ such that $\rho \leq i\leq n-\rho+1$ the following inequality holds {\em
\begin{equation}
\label{L_es_mom1}
\textbf{E}|\xin|^k <C(\rho)\lf \textbf{E}|X_1|^{\delta}g^{-1}\Bigl(\frac i{n+1}\Bigr) \rf^{\rho},
\end{equation} }
where one can put $C(\rho)=2\sqrt{\rho}\exp(\rho+7/6)$.}
 
\bigskip

\noindent{\bf Consequence.} {\it For arbitrary $0<\al<\be<1$ and for all $n$ and $i$ such that $\rho \leq n\al< i <n\be \leq n-\rho +1$, the following inequality holds  {\em
\begin{equation*}
\label{seq}
\textbf{E}|\xin|^k <C(\al,\be,\rho)\lf \textbf{E}|X_1|^{\delta}\rf^{\rho},
\end{equation*} }
where $C(\al,\be,\rho)$ is constant, and depends only on $\al$, $\be$ and $\rho$.
}         

\bigskip
Estimate~\eqref{L_es_mom1} was applied in \citet{gri89,gri94} for obtaining estimates in the Central Limit Theorem for $L$-statistics, and it was applied in \cite{kor_bor93} for the analysis of the asymptotic behavior of random permanent measures for $L$-statistics. It should be noted that many inequalities and other auxiliary results for order statistics are included in \cite{helm_82}, \cite{zwet70}, \cite{shor_weln_86}.

\bigskip
\noindent{\bf Proof of Theorem 1.} To begin with let us note that by Chebyshev's inequality,
$$
\lr F^{-1}(u)\rr^{\delta}(1-u)\le \textbf{E}|X_1|^{\delta}
$$
for $u\ge F(0)$ and
$$
\lv F^{-1}(u)\rv^{\delta}u\le \textbf{E}|X_1|^{\delta}
$$
for $u<F(0)$  where
$$
F^{-1}(u)= \inf \{ x: F(x) \ge u \}.
$$

We first assume that  $\rho+1\le i\le n-\rho$.  Then we obtain
\begin{equation}
\label{L_es_mom2}
\begin{split}
\textbf{E}|\xin|^k =&\frac 1{B(i,n-i+1)}\int_0^1\lv F^{-1}(u)\rv^{k} u^{i-1}(1-u)^{n-i}\, d\, u\\
\le &\frac{\lf \textbf{E}|X_1|^{\delta}\rf^{\rho}
}{B(i,n-i+1)}\lr \int_0^{F(0)} u^{i-\rho-1}(1-u)^{n-i}\, d\, u\rp\\
&+\lp \int_{F(0)}^1 u^{i-1}(1-u)^{n-\rho-i}\, d\, u\rr\\
<&\lf \textbf{E}|X_1|^{\delta}\rf^{\rho}\lr\frac{B(i-\rho,n-i+1)}{B(i,n-i+1)}+ \frac{B(i,n-\rho-i+1)}{B(i,n-i+1)}\rr.
\end{split}
\end{equation}

To evaluate the ratios of beta functions on the right hand side of~\eqref{L_es_mom2}, we apply the inequalities
 \begin{equation}
\label{L_es_mom3}
\sqrt{2\pi}x^{x+1/2}e^{-x}<\Gamma(1+x)<\sqrt{2\pi}x^{x+1/2}e^{-x+1/({12} x)},
\end{equation}
where $x\in\mathbb{R}$, $x>0$, which follows from Stirling's expansion. By using~\eqref{L_es_mom3}, we will prove that for all $\rho>0$ and for all $i\ge \rho+1$
 \begin{equation}
\label{L_es_mom4}
{B(i-\rho,n-i+1)}/{B(i,n-i+1)}<e^{1+7/6}(n/i)^{\rho}.
\end{equation}
Relation~\eqref{L_es_mom4} and the symmetry imply that
 \begin{equation}
\label{L_es_mom5}
{B(i,n-\rho-i+1)}/{B(i,n-i+1)}<e^{1+7/6}(n/(n-i+1))^{\rho}
\end{equation}
for all  $\rho>0$ and for all  $i \le n-\rho$. Together~\eqref{L_es_mom2} and bounds~\eqref{L_es_mom4}-\eqref{L_es_mom5} yield~\eqref{L_es_mom1} in the case $\rho+1\le i\le n-\rho$.
Now we will prove~\eqref{L_es_mom4}. We have
 \begin{equation}
\label{L_es_mom6}
{B(i-\rho,n-i+1)}/{B(i,n-i+1)}=\frac{\Gamma(i-\rho)\cdot\Gamma(n+1)}{\Gamma(i)\cdot\Gamma(n-\rho+1)}.
\end{equation}

First suppose that $i\ge \rho+2$, so therefore $i-\rho-1 \ge 1$, and applying~\eqref{L_es_mom3} on the right hand side of~\eqref{L_es_mom6}, we find that~\eqref{L_es_mom6} is less than
 \begin{equation}
\label{L_es_mom7}
e^{1/6}\lr\frac{(i-\rho-1)n}{(n-\rho)(i-1)} \rr^{1/2}\lr\frac{n}{i-1}\rr^{\rho}\lr\frac{n}{n-\rho}\rr^{n-\rho}\lr\frac{i-\rho-1}{i-1}\rr^{i-\rho-1}.
\end{equation}

To continue let us note that
$$
\lr\frac{n}{n-\rho}\rr^{n-\rho}<e^{\rho},
$$
$$
\lr\frac{(i-\rho-1)n}{(n-\rho)(i-1)} \rr^{1/2}<1, \quad \ \lr\frac{i-\rho-1}{i-1}\rr^{i-\rho-1}<e^{-\rho+\rho^2/(i-1)}, $$
so the magnitude~\eqref{L_es_mom7} is less than
 \begin{equation}
\label{L_es_mom8}
e^{1/6+\rho^2/(i-1)}  \lr\frac{n}{i-1}\rr^{\rho}<e^{\rho+7/6} \lr\frac ni\rr^{\rho}.
\end{equation}

Now we will estimate~\eqref{L_es_mom6} for the case
\begin{equation}
\label{L_es_mom9}
\rho +1\le i\le \rho +2.
\end{equation}

Then we have $1\le i-\rho<2$, $\Gamma(i-\rho)\le 1$, and so the right hand side of~\eqref{L_es_mom6} is not larger than
\begin{equation}
\label{L_es_mom10}
\Gamma(n+1)/(\Gamma(i)\Gamma(n-\rho+1)).
\end{equation}

If $\rho$ is an~integer, then $\Gamma(n+1)/\Gamma(n-\rho+1)<n^{\rho}$ and  $i$ from~\eqref{L_es_mom9} is $i=\rho+1$. For such  $i$ we have
\begin{equation}
\label{L_es_mom11}
\Gamma(n+1)/(\Gamma(i)\Gamma(n-\rho+1))<\frac{n^{\rho}e^{\rho}}{\sqrt{\rho}\rho^{\rho}}<e^{\rho+1}\lr n/i\rr^{\rho},
\end{equation}
and so~\eqref{L_es_mom4} follows. And further, if $\rho$ is not an~integer, then $i$ from~\eqref{L_es_mom9} is  $i=[\rho]+2$ (here and further $[k]$ means the whole part of positive $k$). Suppose at first that $[\rho]=0$. Then  $i=2$ and~\eqref{L_es_mom10} is equal to
\begin{equation}
\label{L_es_mom12}
\Gamma(n+1)/\Gamma(n-\rho+1))< e^{1/12}n^{\rho}\sqrt{1+\rho}< e^{1/12+\rho/2+1}\lr n/i\rr^{\rho},
\end{equation}
and so~\eqref{L_es_mom4} is true. Now we will suppose that $|\rho|\ge 1$, then we have
\begin{equation}
\label{L_es_mom13}
\begin{split}
\Gamma(n+1)/(\Gamma(i)&\cdot\Gamma(n-\rho+1))< e^{1/12}n^{\rho}\sqrt{1+\rho}\, \Gamma^{-1}( [\rho]+2)\\
<& e^{1/12}\sqrt{\frac{1+\rho}{1+[\rho]}}\frac 1{\sqrt{2\pi}}\lr \frac n{1+[\rho]} \rr^{\rho} e^{\rho} \lr\frac{e}{1+[\rho]}\rr^{1+[\rho]-\rho}\\
<& \frac 1{\sqrt{2\pi}} \lr\frac 32\rr^{3/2}e^{\rho+1/12}\lr \frac n {1+[\rho]}\rr^{\rho}\\
<&e^{\rho+1}(n/i)^{\rho}.
\end{split}
\end{equation}

Together bounds~\eqref{L_es_mom6}--\eqref{L_es_mom8} and~\eqref{L_es_mom11}--\eqref{L_es_mom13} imply~\eqref{L_es_mom4}. Thus we have proved our theorem for the case $\rho+1\le i\le n-\rho$ for all $\rho>0$.
It still remains to establish~\eqref{L_es_mom1} for the extreme values of $i$. Now let
\begin{equation}
\label{L_es_mom14}
\rho\le  i < \rho+1.
\end{equation}

If $\rho\le 1$ (i.e. $k\le\delta$), then $i$ from~\eqref{L_es_mom14} is $i=1$. Now by H\"{o}lder's inequality we obtain
\begin{equation*}
\begin{split}
\textbf{E}|X_{1:n}|^k=&B^{-1}(1,n)\int_{-\infty}^{\infty}|x|^k(1-F(x))^{n-1}\, d F(x)\\
\le &B^{-1}(1,n) \lf \int_{-\infty}^{\infty}|x|^{\delta}\, d\, F(x)\rf^{\rho}\\
=&n^{\rho} \lf \textbf{E}|X_1|^{\delta}\rf^{\rho},
\end{split}
\end{equation*}
and so~\eqref{L_es_mom1} is true. Now let $\rho > 1$ (i.e. $k > \delta$). In $\rho$ is an~integer, then $i$ from~\eqref{L_es_mom14} is $i=\rho$. Now for that $i$ we have
\begin{equation}
\label{L_es_mom15}
\begin{split}
\textbf{E}|X_{1:n}|^k=&B^{-1}(\rho,n-\rho+1)\int_{-\infty}^{\infty}|x|^k F^{\rho-1}(x)(1-F(x))^{n-\rho}\, d F(x)\\
\le &B^{-1}(\rho,n-\rho+1)\lf \textbf{E}|X_1|^{\delta}\rf^{\rho -1}  \int_{-\infty}^{\infty}|x|^{\delta}(1-F(x))^{n-2\rho+1}\, d F(x)\\
\le &B^{-1}(\rho,n-\rho+1)\lf \textbf{E}|X_1|^{\delta}\rf^{\rho} .
\end{split}
\end{equation}

Further,
\begin{equation}
\label{L_es_mom16}
\begin{split}
B^{-1}(\rho,n-\rho+1)=&\Gamma(n+1)/(\Gamma(\rho)\cdot\Gamma(n-\rho+1))\\
<& e^{1/12-1} \frac 1{\sqrt{2\pi}}\sqrt{\rho-1}\frac{n^n}{(\rho-1)^{\rho}(n-\rho)^{n-\rho}}\sqrt{\frac{n}{n-\rho}}\\
<& e^{1/12-1}\frac {\sqrt{\rho-1}}{\sqrt{2\pi}} \lr \frac {n}{\rho-1} \rr^{\rho} e^{\rho}\sqrt{2} \\
<& \frac 1{\sqrt{\pi}} e^{1/12+\rho+1}\sqrt{\rho-1} (n/i)^{\rho}.
\end{split}
\end{equation}

Now we assume that  $\rho > 1$  is not an~integer. Then  $i$ from~\eqref{L_es_mom14} is $i=[\rho]+1$. For such $i$ we obtain
\begin{equation}
\label{L_es_mom17}
\begin{split}
\textbf{E}|X_{1:n}|^k=&B^{-1}([\rho]+1,n-[\rho])\int_{-\infty}^{\infty}|x|^k F^{[\rho]}(x)(1-F(x))^{n-[\rho]-1}\, d F(x)\\
\le &\frac{\lf \textbf{E}|X_1|^{\delta}\rf^{[\rho]}}{B([\rho]+1,n-[\rho])}  \int_{-\infty}^{\infty}|x|^{\delta(\rho-[\rho])}(1-F(x))^{n-2[\rho]-1}\, d F(x).
\end{split}
\end{equation}

By H\"{o}lder's inequality the right hand side of~\eqref{L_es_mom17} does not exceed
\begin{equation}
\label{L_es_mom18}
\begin{split}
&\frac{\lf \textbf{E}|X_1|^{\delta}\rf^{[\rho]}}{B([\rho]+1,n-[\rho])}  \cdot
\frac{\lf \textbf{E}|X_1|^{\delta}\rf^{\rho-[\rho]}}{\lf B(1,n-2[\rho])\rf^{\rho-[\rho]-1}}\\
=&\lf \textbf{E}|X_1|^{\delta}\rf^{\rho}\frac{n!(n-2[\rho])^{\rho-[\rho]-1}}{[\rho]!(n-[\rho]-1)!}\\
\le &\lf \textbf{E}|X_1|^{\delta}\rf^{\rho} \frac{n^{\rho}}{[\rho]!}\lr \frac{n}{n-2[\rho]}\rr^{[\rho]+1-\rho}.
\end{split}
\end{equation}

If $[\rho]=1$, then the right hand side of~\eqref{L_es_mom18} is equal to
\begin{equation*}
\lf \textbf{E}|X_1|^{\delta}\rf^{\rho} n^{\rho} \lr 1+ \frac{2}{n-2}\rr^{2-\rho}.
\end{equation*}

Since $n-2\ge 2$ (because  $n\ge 2\rho+1>3$), the latter quantity is less than
\begin{equation*}
\lf \textbf{E}|X_1|^{\delta}\rf^{\rho} n^{\rho}\cdot 2.
\end{equation*}
Thus for  $i=[\rho]+1=2$ we have the inequality
\begin{equation*}
\lf \textbf{E}|X_1|^{\delta}\rf^{\rho} n^{\rho}\cdot 2 <2e^{\rho}\lf\frac {n^2}{2(n-2)}\rf^{\rho},
\end{equation*}
which implies~\eqref{L_es_mom1}. It remains to consider the case $[\rho]\ge 2$. In view of~\eqref{L_es_mom3} the right hand side of~\eqref{L_es_mom18} is less than
\begin{equation*}
\lf \textbf{E}|X_1|^{\delta}\rf^{\rho} \frac 1{\sqrt{2\pi}} e^{[\rho]}\lr\frac{n}{[\rho]}\rr^{\rho}\sqrt{[\rho]}\lr2+\frac
 1{[\rho]}\rr^{[\rho]+1-\rho}<\lf \textbf{E}|X_1|^{\delta}\rf^{\rho}\sqrt{[\rho]}\lr n/[\rho]\rr^{\rho}e^{[\rho]}.
\end{equation*}
So for  $i=[\rho]+1$ we get the estimate
\begin{equation*}
\begin{split}
\textbf{E}|X_{i:n}|^k< &\lf \textbf{E}|X_1|^{\delta}\rf^{\rho}e^{\rho} \sqrt{\rho} \, e^{\frac{\rho-[\rho]}{[\rho]}\,-\,(\rho-[\rho])+1}
\lr\frac ni\rr^{\rho}\\
<&e^{\rho+1}\sqrt{\rho}\lf \textbf{E}|X_1|^{\delta}\rf^{\rho}\lr g^{-1}\Bigl(\frac i{n+1}\Bigr)\rr^{\rho}.
\end{split}
\end{equation*}
Thus, inequality~\eqref{L_es_mom1} is proved for all $i$ such that $\rho\le i<\rho+1$ and all $\rho>0$,
and in a~standard way consisting in the change of the sign of the original random variables, we obtain  that~\eqref{L_es_mom1} is  also valid  for all $i$ such that $n-\rho < i \le n-\rho+1$ and all $\rho>0$.
The theorem is proved.

\bibliographystyle{apalike}
\bibliography{ref_bounds1}
\end{document}